# The Dagum family of isotropic correlation functions

CHRISTIAN BERG[1], JORGE MATEU[2,*] and EMILIO PORCU[1,2,**]

[1]*Department of Mathematics, University of Copenhagen, Universitetsparken 5, DK-2100 Copenhagen, Denmark.* E-mail: *Berg@math.ku.dk*

[2]*Department of Mathematics, Universitat Jaume I, Campus Riu Sec, E-12071 Castellón, Spain.* E-mail: [*]*mateu@mat.uji.es*; [**]*porcu@mat.uji.es*

A function $\rho : [0, \infty) \to (0, 1]$ is a completely monotonic function if and only if $\rho(\|\mathbf{x}\|^2)$ is positive definite on $\mathbb{R}^d$ for all $d$ and thus it represents the correlation function of a weakly stationary and isotropic Gaussian random field. Radial positive definite functions are also of importance as they represent characteristic functions of spherically symmetric probability distributions. In this paper, we analyze the function

$$\rho(\beta, \gamma)(x) = 1 - \left(\frac{x^\beta}{1 + x^\beta}\right)^\gamma, \qquad x \geq 0, \ \beta, \gamma > 0,$$

called the *Dagum function*, and show those ranges for which this function is completely monotonic, that is, positive definite, on any *d*-dimensional Euclidean space. Important relations arise with other families of completely monotonic and logarithmically completely monotonic functions.

## 1. Introduction

Positive definite functions on groups or semigroups have a long history and are present in many applications in probability theory, operator theory, potential theory, moments problems and several other areas. They enter as an important chapter in all treatments of harmonic analysis, and can be traced back to papers by Carathéodory, Herglotz, Bernstein and Matthias, culminating in Bochner's theorem from 1932–1933. See [24] for details.

Schoenberg's theorem ([25]) explains the possibility of constructing rotation-invariant positive definite and conditionally negative definite functions on Euclidean spaces via completely monotonic functions and Bernstein functions. Positive definite radial functions play a very important role in probability theory and statistics. They occur both as the characteristic functions of spherically symmetric probability distributions ([8], and references therein) and as the correlation functions of isotropic random fields on Euclidean spaces ([18]).

Among the valuable contributions to the development of random field theory attributable to eminent scientists, a crucial part is recovered by those regarding the construction of covariance/correlation functions associated with (wide sense) stationary random fields (RF for short). Covariance functions are fundamental in describing the dependence between the values of a spatial RF at locations separated by a lag vector $\mathbf{x} \in \mathbb{R}^d$. The (very often) complex spatial interactions and types of dependencies call for simplifying assumptions and, typically, the spatial phenomenon is considered a realization of a stationary and additionally isotropic RF so that correlations between observations separated by some vector $\mathbf{x} \in \mathbb{R}^d$ become a function of only the Euclidean distance $\|\mathbf{x}\|$.

Among the correlation models proposed in the literature, two are of particular importance for RF theory: the Matérn class ([16]) and the generalized Cauchy class ([12], and references therein).

[27] argues convincingly that the Matérn covariance function possesses desirable mathematical features that make this function very appealing to describe spatial dependencies. As for the generalized Cauchy class, it has been recently revisited by [13], who show heuristically that this class allows for the decoupling of the local and global behavior of the associated weakly stationary and isotropic RF. The Matérn and Cauchy classes are permissible for any dimension of the associated Euclidean space $\mathbb{R}^d$, $d \geq 1$, and only a few other models are available in the literature: among them, it is worth citing the Gaussian model and the exponential model, associated with the so-called Ornstein–Uhlenbeck process or continuous autoregressive of the first order, that are additionally special cases of the Matérn class for particular settings of the parameters indexing this correlation function.

Recently, [23] introduced the Dagum class of radially symmetric correlation functions. Roughly speaking, for $\beta > 0$, $\gamma > 0$, they considered the function

$$\rho(\beta, \gamma)(x) = 1 - \left(\frac{x^\beta}{1 + x^\beta}\right)^\gamma, \qquad x \geq 0, \tag{1.1}$$

and showed permissibility criteria for this function in a three-dimensional Euclidean space. Next, [17] showed that the Dagum class allows the fractal dimension and the Hurst effect of the associated weakly stationary Gaussian RF to be treated independently, by using the procedure suggested in [13]. Comparisons between Dagum and Cauchy models are considered in Section 5.

It is well known that correlation functions must be positive definite, a condition that has been persistently emphasized in the literature. Showing positive definiteness, for a radial function, in any $d$-dimensional Euclidean space, is a nontrivial matter. From the seminal work of Schoenberg ([5, 25, 26]), we know that $\rho(\beta, \gamma)$ is a completely monotonic function if and only if $\rho(\beta, \gamma)(\|\mathbf{x}\|^2)$ is positive definite on $\mathbb{R}^d$ for all $d$. From the point of view of RF theory and modelling, and for the construction of characteristic functions of spherically symmetric distributed random vectors, it is therefore of interest to determine values of the parameters $\beta$, $\gamma$ for which $\rho(\beta, \gamma)$ is completely monotonic.

Our goal has been to characterize the range of parameters $\beta$, $\gamma$ for which $\rho(\beta, \gamma)$ is completely monotonic and this is done in Theorem 9 below. Necessary conditions for complete monotonicity are that $\beta\gamma \leq 1$ and $\beta \leq 2$, while $\beta\gamma \leq 1$ and $\beta \leq 1$ are sufficient. When $\beta\gamma = 1$, the functions are not completely monotonic for $1 < \beta$. However, for $\beta\gamma < 1$, the picture is more complicated: for $1 < \beta < \beta_*$, we find that $\rho(\beta, \gamma)$ is completely monotonic for $0 < \gamma \leq (1-l(\beta))/(\beta+l(\beta))$ for a certain continuous and strictly increasing function $l$ mapping $[1,2]$ onto $[0,2]$ and such that $l(\beta_*) = 1$.

However, we cannot decide if this is the full range of parameters for which $\rho(\beta, \gamma)$ is completely monotonic.

The function $l$, together with another function $c$, is used to give a complete characterization of the parameters $\alpha, \beta \geq 0$ such that

$$f_{\alpha,\beta}(x) = \frac{1}{x^\alpha(1+x^\beta)} \tag{1.2}$$

and all its positive powers are completely monotonic. We believe that this new class of completely monotonic functions may be of some interest in itself.

The paper is organized in the following way. Section 2 contains basic facts about completely monotonic, Stieltjes, logarithmically completely monotonic and Bernstein functions in order to make the exposition more readable. In Section 3, we present the auxiliary family (1.2) and determine the range of parameters for which $f_{\alpha,\beta}$ is completely monotonic. In Section 4 we determine the range of parameters for which the functions (1.2) are logarithmically completely monotonic and we use these criteria to discuss the parameter range for (1.1) to be completely monotonic, hence to obtain permissible correlation functions associated with a weakly stationary Gaussian RF defined on any Euclidean space. In Section 5, we comment on possible comparisons between the Dagum and Cauchy models.

## 2. Families of completely monotonic functions

We recall that a function $f:]0,\infty[\to\mathbb{R}$ is called *completely monotonic* if it is $C^\infty$ and

$$(-1)^n f^{(n)}(x) \geq 0 \qquad \text{for } x > 0, n = 0, 1, \ldots. \tag{2.1}$$

By the celebrated theorem of Bernstein, the set $\mathcal{C}$ of completely monotonic functions coincides with the set of Laplace transforms of positive measures $\mu$ on $[0,\infty[$ (cf. [28]), that is,

$$f(x) = \mathbb{L}\mu(x) = \int_0^\infty e^{-xt}\,d\mu(t), \tag{2.2}$$

where the only requirement on $\mu$ is that $e^{-xt}$ is $\mu$-integrable for any $x > 0$. The set $\mathcal{C}$ is closed under addition, multiplication and pointwise convergence.

In [10], the authors call a function $f:\,]0,\infty[\,\to\,]0,\infty[$ *logarithmically completely monotonic* if it is $C^\infty$ and

$$(-1)^n[\log f(x)]^{(n)} \geq 0 \qquad \text{for } x > 0, n = 1, 2, \ldots, \tag{2.3}$$

which are all the conditions of complete monotonicity of $\log f$, except that $\log f$ need not be non-negative. If we denote the class of logarithmically completely monotonic functions by $\mathcal{L}$, then we have $f \in \mathcal{L}$ if and only if $f$ is a positive $C^\infty$-function such that $-(\log f)' \in \mathcal{C}$.

The functions of class $\mathcal{L}$ have been implicitly studied in [1], and Lemma 2.4(ii) in that paper can be stated as the inclusion $\mathcal{L} \subset \mathcal{C}$, a fact also established in [9].

The class $\mathcal{L}$ can be characterized in the following way, established by Horn [15], Theorem 4.4.

**Theorem 1.** *For a function $f:\,]0,\infty[\,\to\,]0,\infty[$, the following statements are equivalent:*

   (i) $f \in \mathcal{L}$;
   (ii) $f^\alpha \in \mathcal{C}$ for all $\alpha > 0$;
   (iii) $\sqrt[n]{f} \in \mathcal{C}$ for all $n = 1, 2, \ldots$.

From Theorem 1, it follows that $\mathcal{L}$ is closed under products and pointwise convergence since $\mathcal{C}$ is. Note that $x^{-\alpha} \in \mathcal{L}$ for all $\alpha \geq 0$.

We shall now consider an important subclass of $\mathcal{L}$.

**Definition 1.** *A function $f:\,]0,\infty[\,\to\,\mathbb{R}$ is called a Stieltjes function if it is of the form*

$$f(x) = a + \int_0^\infty \frac{\mathrm{d}\mu(t)}{x+t}, \tag{2.4}$$

*where $a \geq 0$ and $\mu$ is a positive measure on $[0,\infty[$. (The conditions on $\mu$ imposed by (2.4) can be summarized as $\int 1/(1+t)\,\mathrm{d}\mu(t) < \infty$.)*

The set $\mathcal{S}$ of Stieltjes functions is a convex cone such that $\mathcal{S} \subset \mathcal{C}$. The inclusion is easy to verify, but we even have $\mathcal{S} \subset \mathcal{L}$ because for a Stieltjes function $f$, we have $f^\alpha \in \mathcal{S} \subset \mathcal{C}$ for any $0 < \alpha \leq 1$ (cf. [4] or the survey lecture [5]).

The formula

$$\frac{1}{x(1+x^2)} = \int_0^\infty \mathrm{e}^{-xt}(1 - \cos t)\,\mathrm{d}t \tag{2.5}$$

shows that $1/x(1+x^2)$ is completely monotonic, but it cannot be a Stieltjes function since it has poles at $\pm \mathrm{i}$ and it follows immediately from (2.4) that a Stieltjes function has a holomorphic extension to the cut plane $\mathbb{C}\setminus\,]-\infty,0]$.

**Definition 2.** *A function $f:\,]0,\infty[\,\to\,[0,\infty[$ is called a Bernstein function if it is $C^\infty$ and $f' \in \mathcal{C}$.*

The set of Bernstein functions is denoted $\mathcal{B}$ and is a convex cone closed under pointwise convergence.

Since a Bernstein function is non-negative and increasing, it has a non-negative limit $f(0+)$. Integrating the Bernstein representation of the completely monotonic function $f'$ gives the integral representation of $f \in \mathcal{B}$

$$f(x) = \alpha x + \beta + \int_0^\infty (1 - e^{-xt}) \, d\nu(t), \tag{2.6}$$

where $\alpha, \beta \geq 0$ and $\nu$, called the *Lévy measure*, is a positive measure on $]0, \infty[$ satisfying

$$\int_0^\infty \frac{t}{1+t} \, d\nu(t) < \infty.$$

The following composition result is useful; see [6] or [5].

**Theorem 2.** *Let $X$ be any of the sets $\mathcal{B}, \mathcal{C}, \mathcal{L}$. Then*

$$f \in X, \ g \in \mathcal{B} \Rightarrow f \circ g \in X.$$

## 3. The auxiliary family F

The main result of this paper, regarding the Dagum family as described in equation (1.1), is obtained by making use of the following family of functions, indexed by the two parameters $\alpha, \beta \geq 0$:

$$\mathrm{F} := \left\{ f_{\alpha,\beta}(x) = \frac{1}{x^\alpha(1+x^\beta)}, x > 0 \right\}. \tag{3.1}$$

We want to discuss the range of parameters for which these functions are completely monotonic.

Since a completely monotonic function has, by (2.2), a holomorphic extension to the half-plane $\Re z > 0$ and since $z^\beta = -1$ in this half-plane if $\beta > 2$, it follows that $\beta \leq 2$ is a necessary condition for $f_{\alpha,\beta} \in \mathcal{C}$.

**Theorem 3.** *The following statements hold:*

(i) $f_{\alpha,\beta} \in \mathcal{C}$ *for* $\alpha \geq 0, 0 \leq \beta \leq 1$;
(ii) $f_{\alpha,\beta} \in \mathcal{C}$ *for* $\alpha \geq \beta/2, 1 < \beta \leq 2$.

**Proof.** Clearly, $x^{-\alpha} \in \mathcal{C}$ for $\alpha \geq 0$. Since $x^\beta \in \mathcal{B}$ for $0 \leq \beta \leq 1$ and the reciprocal of a non-zero Bernstein function is completely monotonic (cf. [6], page 66) we get $1/(1+x^\beta) \in \mathcal{C}$ for these values of $\beta$. Using that $\mathcal{C}$ is closed under products, we see that (i) holds.

Composing the completely monotonic function (2.5) with the Bernstein function $x^\beta, 0 < \beta \leq 1$, we get $f_{\beta, 2\beta} \in \mathcal{C}$ by Theorem 2. Multiplying $f_{\beta, 2\beta}$ by $x^{-\alpha}$, we get (ii). $\square$

We define the function $c\colon [1,2] \to [0,\infty[$ as follows:

$$c(\beta) = \inf\{\alpha \geq 0 \mid f_{\alpha,\beta} \in \mathcal{C}\}. \tag{3.2}$$

By Theorem 3(ii), we get $c(\beta) \leq \beta/2 \leq 1$ and since $\mathcal{C}$ is closed under pointwise convergence, we get $f_{c(\beta),\beta} \in \mathcal{C}$, hence

$$f_{\alpha,\beta} \in \mathcal{C} \Leftrightarrow \alpha \geq c(\beta), \qquad 1 \leq \beta \leq 2. \tag{3.3}$$

While we know from Theorem 3 that $c(1) = 0$ and $c(2) \leq 1$, we next prove that $c(2) = 1$.

**Lemma 1.**

$$f_{\alpha,2} \in \mathcal{C} \Leftrightarrow \alpha \geq 1.$$

**Proof.** We have to show that $f_{\alpha,2} \notin \mathcal{C}$ for $0 \leq \alpha < 1$. This is clear for $\alpha = 0$ since

$$\frac{1}{1+x^2} = \int_0^\infty e^{-xt} \sin t \, dt.$$

For $0 < \alpha$, we define the function $\kappa_\alpha$ by

$$\kappa_\alpha(t) = \frac{t^{\alpha-1}}{\Gamma(\alpha)}, \qquad t > 0, \tag{3.4}$$

and it is easy to see that

$$x^{-\alpha} = \frac{1}{\Gamma(\alpha)} \int_0^\infty e^{-xt} t^{\alpha-1} \, dt,$$

that is, that $x^{-\alpha}$ is the Laplace transform of $\kappa_\alpha$, hence

$$f_{\alpha,2}(x) = \int_0^\infty e^{-xt} \eta_\alpha(t) \, dt, \qquad \eta_\alpha(t) = \frac{1}{\Gamma(\alpha)} \int_0^t (t-s)^{\alpha-1} \sin s \, ds. \tag{3.5}$$

The assertion follows if we prove that $\eta_\alpha(2\pi) < 0$ for $0 < \alpha < 1$. We find

$$\Gamma(\alpha)\eta_a(2\pi) = \int_0^{2\pi} (2\pi - s)^{\alpha-1} \sin s \, ds = \int_0^\pi ((2\pi - s)^{\alpha-1} - (\pi - s)^{\alpha-1}) \sin s \, ds < 0$$

because $(2\pi - s)^{\alpha-1} < (\pi - s)^{\alpha-1}$ and $\sin s > 0$ for $0 < s < \pi$. □

Formula (2.5) can be generalized. For $1 < \beta < 2$, we want to find the function $\varphi_\beta$ such that for $x \geq 0$,

$$\frac{1}{1+x^\beta} = \int_0^\infty e^{-xt} \varphi_\beta(t) \, dt. \tag{3.6}$$

**Theorem 4.** *Let $1 < \beta < 2$. The function $\varphi_\beta$ of (3.6) is given by*

$$\varphi_\beta(t) = \frac{\sin(\beta\pi)}{\pi} \int_0^\infty \frac{e^{-ts} s^\beta \, ds}{1 + 2s^\beta \cos(\beta\pi) + s^{2\beta}}$$
$$- \frac{2}{\beta} \exp\left(t\cos\left(\frac{\pi}{\beta}\right)\right) \cos\left(\frac{\pi}{\beta} + t\sin\left(\frac{\pi}{\beta}\right)\right), \qquad t \geq 0. \qquad (3.7)$$

**Proof.** Fix $t \in \mathbb{R}$. The function $f(z) = e^{tz}/(1 + z^\beta)$ is meromorphic in the cut plane $\mathcal{A} = \mathbb{C} \setminus ]-\infty, 0]$ with simple poles at $z = \exp(\pm i(\pi/\beta))$. Defining

$$\varphi_\beta(t) = \frac{1}{2\pi} \int_{-\infty}^\infty \frac{e^{ity} \, dy}{1 + (iy)^\beta}, \qquad (3.8)$$

we get $\varphi_\beta(t) = 0$ for $t \leq 0$ by integrating $f$ along the half-circle in the right half-plane with diameter $[-iR, iR]$ and letting $R \to \infty$. For $t > 0$, we integrate $f$ along the contour consisting of

$$[-iR, iR], \qquad z = Re^{i\theta}, \qquad \theta \in \left[\frac{\pi}{2}, \pi - \varepsilon\right] \cup \left[-\pi + \varepsilon, -\frac{\pi}{2}\right],$$
$$z = re^{\pm i(\pi - \varepsilon)}, \qquad r \in [0, R],$$

where $R > 1, 0 < \varepsilon < \pi(1 - 1/\beta)$, and next we let $R \to \infty$ and $\varepsilon \to 0$. Since the residues of $f$ at the poles are

$$\mathrm{Res}(f, \exp(\pm i(\pi/\beta))) = -\frac{1}{\beta} e^{\pm i(\pi/\beta)} \exp(te^{\pm i(\pi/\beta)}),$$

we get the expression (3.7) for $\varphi_\beta(t)$. By the Fourier inversion theorem, $1/(1 + (iy)^\beta)$ is the Fourier transform of $\varphi_\beta$ defined by (3.8) and this shows (3.6). $\square$

***Remark 1.*** It is clear that (3.6) holds for $\beta = 1$ with $\varphi_1(t) = \exp(-t)$ and for $\beta = 2$ with $\varphi_2(t) = \sin t$. The latter also follows immediately from (3.7) by letting $\beta \to 2$. However, it is not directly possible to let $\beta \to 1$ in this formula. After an integration by parts using

$$\int \frac{\beta s^{\beta - 1} \, ds}{1 + 2\cos(\beta\pi) s^\beta + s^{2\beta}} = \frac{1}{\sin(\beta\pi)} \arctan \frac{s^\beta + \cos(\beta\pi)}{\sin(\beta\pi)},$$

we find

$$\varphi_\beta(t) = -\frac{1}{\beta\pi} \int_0^\infty \left(\arctan \frac{s^\beta + \cos(\beta\pi)}{\sin(\beta\pi)}\right)(1 - ts)e^{-ts} \, ds$$
$$- \frac{2}{\beta} \exp\left(t\cos\left(\frac{\pi}{\beta}\right)\right) \cos\left(\frac{\pi}{\beta} + t\sin\left(\frac{\pi}{\beta}\right)\right), \qquad t > 0, \qquad (3.9)$$

and splitting the integral into integrals from 0 to 1 and from 1 to $\infty$, it is possible to let $\beta \to 1$ and to obtain $\exp(-t)$ in the limit.

Generalizing (3.5), we have

$$f_{\alpha,\beta}(x) = \int_0^\infty e^{-xt} \eta_{\alpha,\beta}(t) \, dt, \qquad \eta_{\alpha,\beta}(t) = \frac{1}{\Gamma(\alpha)} \int_0^t (t-s)^{\alpha-1} \varphi_\beta(s) \, ds \qquad (3.10)$$

and $f_{\alpha,\beta} \in \mathcal{C}$ if and only if $\eta_{\alpha,\beta}(t) \geq 0$ for $t \in [0, \infty[$.

**Remark 2.** We have $c(\beta) > 0$ for $1 < \beta < 2$. In fact, since $\varphi_\beta(t) < 0$ in the $t$-intervals where $\cos(\frac{\pi}{\beta} + t \sin \frac{\pi}{\beta}) > 0$, we see by continuity that if we fix such a $t$, then also $\eta_{\alpha,\beta}(t) < 0$ for $\alpha$ sufficiently close to 0. Here, we use the fact that $\lim_{\alpha \to 0} \kappa_\alpha = \delta_0$ vaguely (cf. [6], Proposition 14.26).

**Theorem 5.** *The function $c : [1, 2] \to [0, 1]$ is strictly increasing, lower semi-continuous and continuous from the left.*

**Proof.** Fix $1 \leq \beta_1 < \beta_2 \leq 2$. We compose $f_{c(\beta_2), \beta_2} \in \mathcal{C}$ with the function $x^{\beta_1/\beta_2}$ and get

$$f_{c(\beta_2), \beta_2}(x^{\beta_1/\beta_2}) = \frac{1}{x^{\beta_1 c(\beta_2)/\beta_2}(1 + x^{\beta_1})} \in \mathcal{C}$$

by Theorem 2, hence $c(\beta_1) \leq \beta_1 c(\beta_2))/\beta_2 < c(\beta_2)$ as asserted. Fix $1 < \beta_0 \leq 2$ and $0 < \alpha < c(\beta_0)$. Then $f_{\alpha,\beta_0} \notin \mathcal{C}$ and using the fact that $\mathcal{C}$ is closed under pointwise convergence on $]0, \infty[$, we see that $f_{\alpha,\beta} \notin \mathcal{C}$ for $\beta$ in a sufficiently small neighborhood $V$ of $\beta_0$ relative to $]1, 2]$, hence $c(\beta) \geq \alpha$ for $\beta \in V$, but this proves that $c$ is lower semi-continuous at $\beta_0$. An increasing and lower semi-continuous function is continuous from the left. □

**Remark 3.** We conjecture that $c$ is continuous but have not been able to verify this.

## 4. The Dagum family and the class $\mathcal{L}$

In this section, we inspect the range of parameters for which the Dagum class, defined in equation (1.1), is completely monotonic on the positive real line.

Let us begin by showing the range of parameters $\alpha, \beta$ for which $f_{\alpha,\beta} \in \mathcal{L}$.

**Theorem 6.** *The following statements hold:*

(i) $f_{\alpha,\beta} \in \mathcal{L}$ *for* $\alpha \geq 0, 0 \leq \beta \leq 1$;
(ii) $f_{\alpha,2} \in \mathcal{L} \Leftrightarrow \alpha \geq 2$.

**Proof.** (i) By Theorem 1(iii), it suffices to prove that $(x^\alpha(1 + x^\beta))^{-1/n} \in \mathcal{C}$ for $n = 1, 2, \ldots$, but this holds because $(1 + x^\beta)^{1/n} \in \mathcal{B}$ by Theorem 2 and the reciprocal of a Bernstein function is completely monotonic.

(ii) We shall determine the parameters $\alpha \geq 0$ such that $-(\log f_{\alpha,2})' \in \mathcal{C}$, that is, such that

$$\frac{\alpha}{x} + \frac{2x}{1+x^2} \in \mathcal{C}.$$

However, this function is indeed the Laplace transform of $\alpha + 2\cos t$, which is non-negative on $[0, \infty[$ if and only if $\alpha \geq 2$. $\square$

We now define the function $l : [1, 2] \to [0, \infty[$ by

$$l(\beta) = \inf\{\alpha \geq 0 \mid f_{\alpha,\beta} \in \mathcal{L}\}. \tag{4.1}$$

Using the fact that $\mathcal{L}$ is closed under pointwise convergence, we get $f_{l(\beta),\beta} \in \mathcal{L}$ and observing that if $f_{\alpha,\beta} \in \mathcal{L}$, then $f_{\alpha',\beta} \in \mathcal{L}$ for any $\alpha' \geq \alpha$, we get that

$$f_{\alpha,\beta} \in \mathcal{L} \Leftrightarrow \alpha \geq l(\beta), \qquad 1 \leq \beta \leq 2. \tag{4.2}$$

From Theorem 6, we know that $l(1) = 0$, $l(2) = 2$ and, clearly, $c(\beta) \leq l(\beta)$.
For $1 \leq \beta \leq 2$, we introduce the function

$$\psi_\beta(t) = \int_0^t \varphi_\beta(s)\, ds, \qquad t \geq 0. \tag{4.3}$$

Since it is equal to $\kappa_1 * \varphi_\beta$, we know that

$$\frac{1}{x(1+x^\beta)} = \mathbb{L}[\psi_\beta](x), \tag{4.4}$$

which is completely monotonic because $c(\beta) \leq 1$, hence $\psi_\beta(t) \geq 0$ for $t \geq 0$. From formula (3.7), we find $(1 < \beta < 2)$

$$\psi_\beta(t) = \rho_\beta(t) + \tau_\beta(t), \tag{4.5}$$

with

$$\rho_\beta(t) = 1 - \frac{2}{\beta} \exp\left(t \cos \frac{\pi}{\beta}\right) \cos\left(t \sin \frac{\pi}{\beta}\right) \tag{4.6}$$

and

$$\tau_\beta(t) = -\frac{\sin(\beta\pi)}{\pi} \int_0^\infty \frac{e^{-ts} s^{\beta-1}\, ds}{1 + 2s^\beta \cos(\beta\pi) + s^{2\beta}} \tag{4.7}$$

because $\psi'_\beta = \varphi_\beta$ and $\psi_\beta(0) = 0$. Note that $\tau_\beta(t) \in \mathcal{C}$ with $\tau_\beta(0) = 2/\beta - 1$ and

$$\rho_\beta\left(\frac{\pi}{\sin \pi/\beta}\right) = 1 + \frac{2}{\beta} \exp\left(\pi \cot \frac{\pi}{\beta}\right) > 1,$$

hence $\psi_\beta(t)$ has a global maximum which is greater than 1 on $[0, \infty[$.

Substituting $s^\beta = -u\sin(\beta\pi) - \cos(\beta\pi)$ in (4.7), we get

$$\tau_\beta(t) = \frac{1}{\beta\pi} \int_{-\cot(\beta\pi)}^{\infty} \frac{\exp(-t(-u\sin(\beta\pi) - \cos(\beta\pi))^{1/\beta})}{1+u^2}\,du. \tag{4.8}$$

Let $U = UCB([0,\infty[)$ denote the Banach space of uniformly continuous and bounded real-valued functions on $[0,\infty[$ with the uniform norm. For each $1 < \beta < 2$, it is clear that $\psi_\beta \in U$ since

$$\lim_{t\to\infty} \psi_\beta(t) = \int_0^\infty \varphi_\beta(s)\,ds = 1.$$

**Theorem 7.** *The function $\beta \to \psi_\beta$ from $]1,2[$ to $U$ is continuous. Furthermore,*

$$\lim_{\beta\to 1} \psi_\beta(t) = \psi_1(t) = 1 - e^{-t} \qquad \text{in } U \tag{4.9}$$

*and*

$$\lim_{\beta\to 2} \psi_\beta(t) = \psi_2(t) = 1 - \cos t \qquad \text{locally uniformly for } t \in [0,\infty[. \tag{4.10}$$

**Proof.** The function $\rho_\beta(t)$ is well defined and continuous for $(\beta,t) \in [1,2] \times [0,\infty[$ with

$$\rho_1(t) = 1 - 2e^{-t}, \qquad \rho_2(t) = 1 - \cos t.$$

It is easy to see that $\beta \to \rho_\beta$ is continuous on $[1,2[$ with values in the Banach space $U$, while $\rho_\beta(t) \to \rho_2(t)$ locally uniformly on $[0,\infty[$ when $\beta \to 2$. (We cannot have global uniform convergence since $\rho_\beta$ tends to 1 at infinity while $\rho_2$ is oscillating.) By (4.8), it is easy to see that $\beta \to \tau_\beta$ is continuous from $]1,2[$ to $U$ and that

$$\lim_{\beta\to 1} \tau_\beta(t) = e^{-t}, \qquad \lim_{\beta\to 2} \tau_\beta(t) = 0 \qquad \text{in } U,$$

that is, uniformly for $t \in [0,\infty[$. $\square$

**Corollary 1.** *The function $\beta \to \Psi(\beta) = \max_{t\geq 0} \psi_\beta(t)$ is continuous on $[1,2]$.*

**Proof.** We know that $\Psi(1) = 1$ and that $\Psi(\beta) > 1$ for $1 < \beta \leq 2$. Furthermore,

$$\Psi(\beta) \leq \max_{t\geq 0} \rho_\beta(t) + \max_{t\geq 0} \tau_\beta(t) \leq 1 + \frac{2}{\beta} + \frac{2}{\beta} - 1 = \frac{4}{\beta}.$$

If $F: X \to U$ is a continuous mapping of a topological space $X$ with values in the Banach space $U$, then $\sup F: X \to \mathbb{R}$ is also continuous. Therefore, $\Psi$ is continuous on $[1,2[$. By (4.10), we get

$$\max_{0\leq t\leq 2\pi} \psi_\beta(t) \to \max_{0\leq t\leq 2\pi} \psi_2(t) = 2,$$

hence $\liminf_{\beta \to 2} \Psi(\beta) \geq 2$. On the other hand, $\limsup_{\beta \to 2} \Psi(\beta) \leq \lim_{\beta \to 2} \frac{4}{\beta} = 2$ and a combination of these results gives $\lim_{\beta \to 2} \Psi(\beta) = 2$. □

**Theorem 8.** *The function $l$ is given by*

$$l(\beta) = \beta(\Psi(\beta) - 1), \qquad 1 \leq \beta \leq 2, \tag{4.11}$$

*and it is continuous and strictly increasing. Furthermore, $\Psi$ is increasing.*

**Proof.** We have

$$-(\log f_{\alpha,\beta})'(x) = \frac{\alpha + \beta}{x} - \frac{\beta}{x(1+x^\beta)} = \mathbb{L}[(\alpha+\beta) - \beta\psi_\beta](x)$$

and therefore $f_{\alpha,\beta} \in \mathcal{L}$ if and only if $\psi_\beta(t) \leq 1 + \alpha/\beta$ for $t \geq 0$, which holds if and only if

$$\alpha \geq \beta(\Psi(\beta) - 1),$$

giving the formula (4.11). The continuity of $l$ follows from Corollary 1.

We prove that $l(\beta)/\beta$ is increasing, which implies that $l(\beta)$ is strictly increasing and that $\Psi$ is increasing. Let $1 \leq \beta_1 < \beta_2 \leq 2$. We compose $f_{l(\beta_2),\beta_2} \in \mathcal{L}$ with the function $x^{\beta_1/\beta_2}$ and get

$$f_{l(\beta_2),\beta_2}(x^{\beta_1/\beta_2}) = \frac{1}{x^{\beta_1 l(\beta_2)/\beta_2}(1 + x^{\beta_1})} \in \mathcal{L}$$

by Theorem 2, hence $l(\beta_1) \leq \beta_1 l(\beta_2)/\beta_2 < l(\beta_2)$ as asserted. □

Using that the function $l$ is continuous and strictly increasing such that $l(1) = 0, l(2) = 2$, the following definition makes sense.

***Definition 3.*** *Let $\beta_*$ denote the unique number in $]1,2[$ such that $l(\beta_*) = 1$.*

The number $\beta_*$ is where the graphs of $\Psi(\beta)$ and $1 + 1/\beta$ intersect, and from Figure 1 (produced using Mathematica), we estimate that $\beta_*$ is approximately 1.74. For $\beta_* \leq \beta \leq 2$, we have $l(\beta) > c(\beta)$. It is reasonable to believe that this inequality holds also for $1 < \beta < \beta_*$. In any case, for $1 < \beta \leq 2$, we have

$$f_{\alpha,\beta} \in \mathcal{C} \setminus \mathcal{L} \qquad \text{for } c(\beta) \leq \alpha < l(\beta) \tag{4.12}$$

and

$$f_{\alpha,\beta}^\gamma \in \mathcal{C} \qquad \text{for } \gamma > 0, \alpha \geq l(\beta). \tag{4.13}$$

We are now able to state the main result of this paper, regarding the Dagum family $\rho(\beta,\gamma)$ as defined in (1.1).

**Theorem 9.** *Let $\beta, \gamma > 0$. We then have*

$$\rho(\beta,\gamma) \in \mathcal{C} \Leftrightarrow \frac{x^{\beta\gamma-1}}{(1+x^\beta)^{\gamma+1}} \in \mathcal{C}. \tag{4.14}$$

*If these conditions hold, then $\beta\gamma \leq 1$ and $\beta \leq 2$. In the opposite direction, we have*

(i) *if $\beta\gamma \leq 1$ and $\beta \leq 1$ then $\rho(\beta,\gamma) \in \mathcal{C}$;*
(ii) *if $\beta\gamma = 1$, then*

$$\rho(\beta, 1/\beta) \in \mathcal{C} \Leftrightarrow \beta \leq 1; \tag{4.15}$$

(iii) *if $\beta\gamma < 1$ and $1 < \beta < \beta_*$, then $\rho(\beta,\gamma) \in \mathcal{C}$ if, in addition,*

$$0 < \gamma \leq \frac{1 - l(\beta)}{\beta + l(\beta)}. \tag{4.16}$$

**Proof.** Using the fact that $\rho(\beta,\gamma)$ is a positive $C^\infty$-function, it is completely monotonic if and only if $-\rho(\beta,\gamma)' \in \mathcal{C}$, and a small calculation leads to the right-hand side of (4.14).

If this function is completely monotonic, the value at $x = 0$ has to be strictly positive or infinite and hence $\beta\gamma \leq 1$. Furthermore, $\beta \leq 2$, otherwise the function has poles in the right half-plane and then it cannot be completely monotonic.

If $\beta\gamma \leq 1$ and $\beta \leq 1$, then

$$\frac{x^{\beta\gamma-1}}{(1+x^\beta)^{\gamma+1}} = \left(\frac{1}{x^{(1-\beta\gamma)/(1+\gamma)}(1+x^\beta)}\right)^{1+\gamma} \in \mathcal{C} \tag{4.17}$$

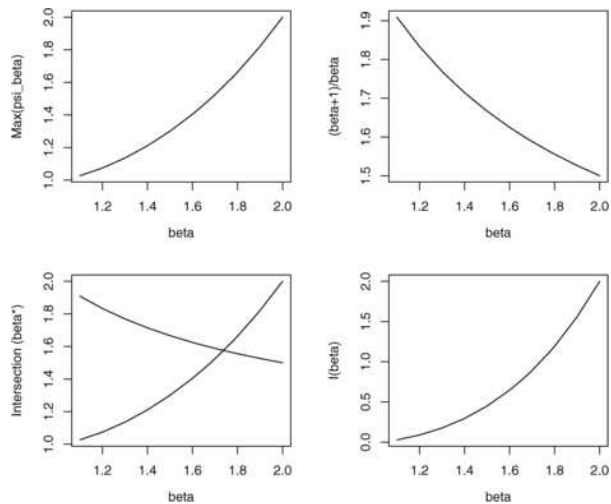

**Figure 1.** From upper left to lower right: $\max(\Psi(\beta))$ vs. $\beta$; $1 + 1/\beta$ vs. $\beta$; intersection of the first two graphs; $l(\beta)$ vs. $\beta$.

because $f_{\alpha,\beta} \in \mathcal{L}$ for any $\alpha \geq 0$, which shows (i).

To see (ii), it suffices, for $c > 0$, to prove the following implication:

$$f_{0,\beta}^c \in \mathcal{C} \Rightarrow \beta \leq 1. \tag{4.18}$$

Assuming $f_{0,\beta}^c \in \mathcal{C}$, we know that $(f_{0,\beta}^c)'' \geq 0$, hence

$$\beta c t^{\beta-2}(1+x^\beta)^{-c-2}[x^b(\beta c + 1) + 1 - \beta] \geq 0, \qquad x > 0,$$

and therefore

$$\lim_{x \to 0+}[x^\beta(\beta c + 1) + 1 - \beta] \geq 0.$$

This shows that $\beta \leq 1$.

The assertion (iii) also follows from formula (4.17) because $l(\beta) < 1$ for $1 < \beta < \beta_*$ and

$$l(\beta) \leq \frac{1-\beta\gamma}{1+\gamma} \Leftrightarrow \gamma \leq \frac{1-l(\beta)}{\beta+l(\beta)}. \qquad \square$$

Theorem 9 shows the range of completely monotonic permissibility for the Dagum class. Under these conditions, we have that this function is the correlation function of an isotropic and Gaussian RF. Rephrasing this, the result shows that there exists a (wide sense) stationary Gaussian RF defined on any $d$-dimensional Euclidean space, whose correlation function is defined by the Dagum class. The same result admits an interpretation in terms of characteristic functions of positive spherically symmetric distributed random vectors of $\mathbb{R}^d$ for all $d \in \mathbb{Z}_+$.

*Remark 4.* For $\beta = 2$, results about the complete monotonicity of

$$g_{\alpha,\lambda}(x) = \frac{1}{x^\alpha(1+x^2)^\lambda}$$

have been obtained in a number of papers. We mention the following results:

(i) $g_{\alpha,\lambda} \in \mathcal{C}$ for $\alpha \geq \lambda \geq 1$ (cf. [11]);
(ii) $g_{\alpha,\lambda} \in \mathcal{C}$ for $\alpha \geq 2\lambda \geq 0$ (cf. [2, 3]);
(iii) $g_{1,\lambda} \in \mathcal{C}$ for $0 \leq \lambda \leq 1$ (cf. [19]);
(iv) $g_{\alpha,\lambda} \notin \mathcal{C}$ for $0 \leq \alpha < \lambda$ (cf. [19, 29]) – the simple proof in [29], Lemma 8, carries over to this situation;
(v) $g_{\alpha,\alpha} \notin \mathcal{C}$ for $0 < \alpha < 1$ (cf. [19]).

## 5. Comparisons between models and some conclusions

In this paper, we have presented some results for the range of parameters allowing the complete monotonicity of the Dagum function. Some comments are in order.

A first consideration is that there is not, to the knowledge of the authors, any comparative study available in the literature regarding the comparison of the Cauchy and the Matérn models. The reason for this is that the Matérn class possesses an explicit analytical form for the associated spectral density and this allows some important properties of the associated random field (RF) to be known. Such properties cannot be established in the case of the Cauchy model, as the Fourier pair can be calculated only for some particular parameter setting. The same holds for the Dagum class. In addition, the Matérn and Cauchy correlation functions represent the most widely used families for geostatistical modelling.

A natural comparison arises between the Dagum and Cauchy correlations, being both models *decouplers* of the fractal dimension and Hurst effect. With this purpose, in a recent paper ([17]) an extensive simulation study is performed in order to highlight the attitude of the Dagum correlation function with respect to decoupling. In this case, the authors use the convenient parametrization

$$\rho(t) = 1 - \left(\frac{t^\gamma}{1+t^\gamma}\right)^{\varepsilon/\gamma}, \qquad t := \|\mathbf{h}\| \in \mathbb{R}^d.$$

One can verify that sufficient conditions for the permissibility of this function on any $d$-dimensional Euclidean space are $\gamma \in (0,2]$ and $\varepsilon < \gamma$. Recall that the Cauchy model has expression

$$\rho(t) = (1+t^\theta)^{-\eta/\theta}, \qquad t := \|\mathbf{h}\| \in \mathbb{R}^d, \ \theta \in (0,2], \eta > 0.$$

Observe that the parameters identifying the fractal dimensions associated, respectively, with Dagum and Cauchy classes, are $\varepsilon, \theta \in (0,1]$, while those identifying the Hurst effect are (respectively) $\delta, \eta \in (0,1)$. From the results contained therein, it is quite evident that, when simulating profiles (i.e., in the one-dimensional case) the Cauchy class shows better performances in comparison with the Dagum class in terms of decoupling.

The dilemma arises in two- or higher-dimensional spaces, where it is far from clear which of the two models performs better. Reaching conclusions, in this case, is not a trivial matter, as several variables can influence the performance of two decouplers in two- or higher-dimensional spaces and the nature of these variables is very often of a non-statistical nature. This is verified, for instance, in [20, 21]. Thus, everything seems to converge to some open problems whose solution could be found through heuristic arguments and deeper simulation studies, taking into account the influence of some real variables on the behavior of the decouplers in bidimensional spaces.

In addition, and as far as estimation of the parameters is concerned, we performed several simulation studies using likelihood approximations such as the *weighted composite likelihood*, proposed by [7], and for both models the results are very promising in the sense that, in both cases, the parameters identifying fractal dimension and Hurst effect are estimated with a high level of accuracy.

Thus, the Dagum model behaves similarly to the Cauchy model in terms of decoupling and estimation properties. In addition, the Dagum function is actually opening several avenues of research in both applied and theoretical studies. For the latter, initial evidence of this fact can be found in [14], where the Cauchy and Dagum functions are studied in terms of extremal rays of positive definite functions. As an example of the former, the Dagum model is lately being considered in material science and mechanical engineering (see, e.g., [21, 22]), where practical analysis in this field is performed using the Dagum function. These considerations reinforce the added value of the Dagum model, which can be considered a good alternative to such other important published models.

## Acknowledgements

This work was partially funded by Grant MTM2007-62923 from the Spanish Ministry of Science and Education. The referee and the Editors are acknowledged with thanks for their thorough reviews that have allowed an earlier version of the manuscript to be considerably improved.